# IDENTIFICATION OF REPETITIVE PROCESSES AT STEADY- AND UNSTEADY-STATE: TRANSFER FUNCTION

Ricardo Antunes[1], Vicente A. González[2], and Kenneth Walsh[3]

## ABSTRACT

Projects are finite terminating endeavors with distinctive outcomes, usually, occurring under transient conditions. Nevertheless, most estimation, planning, and scheduling approaches overlook the dynamics of project-based systems in construction. These approaches underestimate the influence of process repetitiveness, the variation of learning curves and the conservation of processes' properties. So far, estimation and modeling approaches have enabled a comprehensive understanding of repetitive processes in projects at steady-state. However, there has been little research to understand and develop an integrated and explicit representation of the dynamics of these processes in either transient, steady or unsteady conditions. This study evaluates the transfer function in its capability of simultaneously identifying and representing the production behavior of repetitive processes in different state conditions. The sample data for this research comes from the construction of an offshore oil well and describes the performance of a particular process by considering the inputs necessary to produce the outputs. The result is a concise mathematical model that satisfactorily reproduces the process' behavior. Identifying suitable modeling methods, which accurately represent the dynamic conditions of production in repetitive processes, may provide more robust means to plan and control construction projects based on a mathematically driven production theory.

## KEYWORDS

Production, process, system identification, transfer function, system model, theory;

## INTRODUCTION

Construction management practices often lack the appropriate level of ability to handle uncertainty and complexity (Abdelhamid, 2004; McCray and Purvis, 2002) involved in project-based systems resulting in projects failures in terms of projects schedule and budget performance, among other measures (Mills, 2001). Traditional scheduling approaches in construction, such as critical path method, have been unrestrictedly used producing unfinished and erratic plans (Abdelhamid, 2004, Bertelsen, 2003a) consequently creating distrust, and often being abandoned by those conducting project work. Even more recent scheduling approaches, such as the ones based on the line-of-balance method, assume that the production in construction operates at steady-state with constant production rates (Arditi et al. 2001; Lumsden, 1968), where any deviation is understood as variability (Poshdar et al. 2014). However, "the assumption that production rates of construction projects and processes are linear may be erroneous" (Lutz and Hijazi, 1993). Production throughput is highly variable in construction projects (Gonzalez et al. 2009), has transients (Lutz and Hijazi, 1993), occasionally is at unsteady-state (Bernold, 1989, Walsh et al. 2007) and frequently is nonlinear (Bertelsen, 2003b). As such, approaches that depend on constant production rates, *i.e.*, a steady system, possibly produce erroneous and imprecise outcomes.

[1] Ph.D. candidate, Department of Civil and Environmental Engineering - University of Auckland, New Zealand, rsan640@aucklanduni.ac.nz

[2] Senior Lecturer, Department of Civil and Environmental Engineering - University of Auckland, New Zealand, v.gonzalez@auckland.ac.nz

[3] Dean, SDSU-Georgia, San Diego State University, Tbilisi, Georgia, kwalsh@mail.sdsu.edu

The dynamics of the production system in construction is frequently overlooked (Bertelsen, 2003b), and the transient phase is ignored (Lutz and Hijazi, 1993). The general construction management makes no distinction between the production dynamics and disturbance, considering both as variability (Poshdar et al. 2014). However, dynamics, disturbance, and variability have different meanings and action approaches. The dynamics is an essential characteristic of any process, representing the effects of the interaction of components in a system. Process dynamics should be understood, managed and optimized. External factors cause disturbance, which must be filtered, mitigated and avoided consequently reducing any impact on the process, *e.g.*, risk management (Antunes and Gonzalez, 2015). The understanding of these concepts is fundamental to the development of mathematical relations and laws suitable to the construction production system. At this time, construction adopts the manufacturing model, dismissing the application of mathematical approaches to model and manage its production system (Bertelsen, 2003a, Laufer, 1997, McCray and Purvis, 2002).

Although much work has been done to date on production estimates of repetitive processes, more studies need to be conducted to understand and develop the dynamics of these processes. The purpose of this study is to evaluate the transfer function in its capability of identifying and describing the dynamics of project-driven systems in repetitive processes in construction. This topic was identified as being of importance to point out a unique mathematical representation of project-based systems process in transient, unsteady-, and steady-state, furthermore, overcoming a major limitation of fixed production rates estimation approaches. The understanding of project dynamics should improve estimation accuracy approaches and support suitable derivations of manufacturing management practices in order to increase productivity in construction projects. This study is a step towards the development of a mathematically driven production theory for construction.

## A SYSTEM VIEW

Mathematical models have enabled a comprehensive understanding of production mechanisms supporting practices to improve production in manufacturing. Hopp and Spearman (1996) committed to the comprehension of the manufacturing production system. The system approach or system analysis was the problem-solving methodology of choice (Hopp and Spearman, 1996). The first step of this methodology is a system view. In the system view, the problem is observed as a system established by a set of subsystems that interact with each other. Using the system approach, Hopp and Spearman elaborated significant laws to queue systems and the general production in manufacturing. The conservation of material (Wallace J Hopp; Mark L Spearman; Richard Hercher 1996) and capacity laws (Hopp and Spearman, 1996) are particularly attractive, not only according to their importance, but also because they explicitly state one or more system restrictions. These laws place reliance on stable systems, with long runs and at steady-state conditions. However, production in project-based systems, such as construction, involves a mix of processes in steady- and unsteady-state, short and long production runs, and different learning curves (Antunes and Gonzalez, 2015). Hence, unless a construction process fulfills the stability and steady-state conditions, the manufacturing model and, consequently, the laws do not accurately represent production in construction. Alternatively, variants of manufacturing laws must be developed to production in project-based systems that not fulfill those requirements. In this scenario of variety, it is crucial distinguishing between project-based systems conditions, comprehending process dynamics and its behavior.

## SYSTEM IDENTIFICATION

The objective of system identification is to build mathematical models of dynamic systems using measured data from a system (Ljung, 1998). There are several system identification approaches to model different systems, for instance, transfer function. The transfer function is particularly useful because it provides an algebraic description of a system as well means to calculate parameters of the system dynamics and stability. Nevertheless, the modeling capability of the transfer function in construction must be evaluated and tested. In this study, the modeling approach, *i.e.*, transfer function, focuses on replicating the input/output "mapping" observed in a sample data. When the primary goal is the most accurate replication of data, regardless of the mathematical model structure, a black-box modeling approach is useful. Additionally, black-box modeling supports a variety of models (Bapat, 2011; Billings, 2013), which have traditionally been practical for representing dynamic systems. It means that at the end of the modeling, a mathematical description represents the actual process performance rather than a structure biased by assumptions and restrictions. Black-box modeling is a trial-and-error method, where parameters of various models are estimated, and the output from those models is compared to the results with the opportunity for further refinement. The resulting models vary in complexity depending on the flexibility needed to account for both the dynamics and any disturbance in the data. The transfer function is used in order to show the system dynamics explicitly.

### TRANSFER FUNCTION

The transfer function of a system, $G$, is a transformation from an input function into an output function, capable of describing an output (or multiple outputs) by an input (or multiple inputs) change, $y(t) = G(t) \star u(t)$. Although generic, the application of the transfer function concept is restricted to systems that are represented by ordinary differential equations (Mandal, 2006). Ordinary differential equations can represent most dynamic systems in its entirety or at least in determined operational regions producing accurate results (Altmann and Macdonald, 2005; Mandal, 2006). As a consequence, the transfer function modeling is extensively applied in the analysis and design of systems (Ogata, 2010). A generic transfer function makes possible representing the system dynamics by algebraic equations in the frequency domain, $s$. In the frequency domain, the convolution operation transforms into an algebraic multiplication in $s$, which is simpler to manipulate. Mathematically, "the transfer function of a linear system is defined as the ratio of the Laplace transform of the output, $y(t)$, to the Laplace transform of the input, $u(t)$, under the assumption that all initial conditions are

795 Proceedings IGLC-23, July 2015 | Perth, Australia

zero" (Mandal, 2006), Equation 1. Where the highest power of *s* in the denominator of the transfer function is equal to *n*, the system is called a *n*th-order system.

$$G(s) = \frac{\mathcal{L}[u(t)]}{\mathcal{L}[y(t)]} = \frac{U(s)}{Y(s)}$$

*Equation 1: Transfer function*

### TRANSIENT STATE, STEADY-STATE, AND UNSTEADY-STATE RESPONSE

Two parts compose a system response in the time domain, transient, and steady- or unsteady-state. Transient is the immediate system response to an input from an equilibrium state. After the transient state, a system response can assume a steady- or unsteady-state. In a stable system, the output tends to a constant value when $t\to\infty$ (Mandal, 2006). When the system response enters and stays in the threshold around the constant value the system reached the steady-state (Mandal, 2006). The time the stable system takes to reach the steady-state is the settling time, $t_s$. On the other hand, if the response never reaches a final value or oscillates surpassing the threshold when $t\to\infty$ the system is then at unsteady-state. Consequently, the system outputs at unsteady-state vary with time during the on-time interval even induced by an invariable input.

### METHODOLOGY

A sample of 395 meters of continuous drilling was randomly selected from the project of an offshore oil well construction, constituting the process to be modeled. The information containing the drill ahead goal and the current process duration was collected from operational reports and resampled to 181 samples representing the hourly process behavior when commanded by the input, establishing a system. Next, the estimation of a transfer function was used for the determination of a model that represents the dynamics of the system-based process. The estimation uses nine partitions of the dataset creating models based on different data sizes. The best model from each of the nine partitions presenting the lowest estimation unfitness value were selected and cross-validated by the remaining data. Later, the system response of the best model was analyzed.

## CASE STUDY: DRILLING AS A SYSTEM

The subject of this study is the drilling process on a particular offshore well construction project in Brazilian pre-salt. This process was chosen given its high level of repetitiveness. The vertical dimension of repetitiveness is the repetition of the process in the project, *i.e.*, the drilling occurs more than one time in the construction of a well. The horizontal dimension is the repetition of the process in different projects, *i.e.*, the drilling occurs on every well construction project. Such degree of repetitiveness eases comparison and data validation because a repetitive process tends to present patterns in smaller data portions. The case documentation provides details about inputs, outputs and brief explanations of the process parameters. Nonetheless, the documentation does not include any mathematical representation of the processes other than the drilling parameters and other activities performed while drilling, which constitute subsystems. For instance, the work instructions to drill a segment of 28 meters on seabed:

- Drill ahead 8 1/2" hole from 3684 m to 3712 m with 480 gpm, 1850 psi, 15-25k WOB, 120 rpm, 15-20 kft.lbs torq. Perform surveys and downlinks as per directional driller instructions. Pump 15 bbl fine pill and 50 bbl hi-vis pill every two stand as per mud engineer instructions.

The primary input, 'drill ahead from 3684 m to 3712 m', and the parameters, such as torque and rpm, directly affect the drilling process. However, the system view unifies the different parts of the system, *i.e.*, the subsystems, into an effectual unit using a holistic perspective. The holistic perspective allows the creation of a system driven by a primary input while all others variables interact as subsystems of the main system.

To fully establish a system, an input has to be applied to the process in order produce an output. Accordingly, the input, $u(t)$, consists of a drilling ahead depth goal, *e.g.*, 3712 meters, that is applied to the drilling process, producing the output, $y(t)$, that is the actual depth, in meters. In the process, $G(t)$, the drilling crew responds to the drilling ahead goal by drilling and performing related tasks, which increases the actual well depth over time until reaching the drilling depth goal. Then, a new drill ahead goal is set, and the process performs the cycle. The sample data corresponds with a well depth increase from 3305 to 3700 meters at a variable rate based on operational choices. The 181 samples represent the hourly input and process behavior response, *i.e.*, output. The input and output data are cumulative due to physical restrictions. In other words, it is impossible to drill from 3435 meters without prior drilling from the seabed at 924 meters from water level to 3435 meters in the hole. Consequently the drill ahead goal as well as the actual depth values are always greater than the previous values. Figure 1 displays the general system representation of the process with its measured input and output, drill ahead goal and actual depth respectively. Two criteria guided the choice of drilling goal as input. The first criterion is that the drill ahead goal is the primarily directive to achieve the objective of the project, setting the peace to build the well. The drill ahead goal is an adaptive plan in which the team has to examine the current conditions of the well and determine the best drill ahead goal. It relies on guidelines and procedures but in the end its a human decision. The second criterion is that this particular arrangement illustrates the number of items arriving in a queuing system at time *t*. Additionally, the output represents the number of items departing in the queuing system at time *t*. Such input-output arrangement is instrumental to adapt manufacturing-based models such as the Little's Law (Little, 1961) to construction in further research.

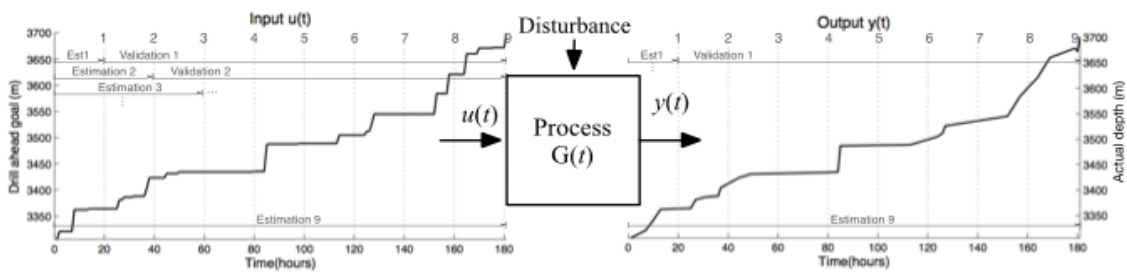

*Figure 1: System representation of the case study*

### INITIAL MODELING APPROACH

A simple model is attempted initially before progressing to more complex structures until reaching the required model accuracy. Simpler models are easier to interpret, a desired feature in this study. However, if that model unsatisfactorily simulates the measured data, it may be necessary to use more complex models. The simpler system identification approach is the transfer function. Hence, transfer function might be a good starting point in order to identify, model and understand the behavior of a system. The sample data was partitioned in nine combinations representing nine stages in time, as shown in Figure 1. The first partition is at the 20th hour. Therefore, the data from zero to 20 hours was used as estimation data for $G_1$, and from 21 to 181 hours as the validation data. The second partition happens at 40 hours mark. In the same way, the estimation for $G_2$ is composed of the data from zero to 40 hours mark, and the validation is from the 41st to the 181st hour. This pattern repeats until the 180-time stamp. At this partition, almost the whole sample constitutes the estimation data, and only one sample is left for validation of $G_9$. The model from this partition, $G_9$, merely fits the estimation data once there is virtually no data that could be used to validate the model. Based on black-box trial-and-error approach, the model parameters of the transfer function of first-order (Ogata, 2010) were generated for each partition using the iterative prediction-error minimization algorithm (Ljung, 2010) from MATLAB's System Identification Toolbox. A first-order transfer function eases the model interpretability.

## MODEL DEVELOPMENT

Three transfer functions, which showed the lowest unfitness values, calculated by 100% – normalized root mean square (NRMSE) (Armstrong and Collopy, 1992; Ljung, 2010), were selected for each of the nine data partitions, constituting the best models. A perfect fit corresponds to zero meaning that the simulated or predicted model output is exactly the same as the measured data.

## MODEL QUALITY ASSESSMENT

The initial models were later refined using the prediction-error minimization algorithm (Ljung, 1998). After refinement, the models that achieved the lowest unfitness values to each estimation data partition that they derived from were then validated using the remaining data of their partition. Figure 2 shows the quality measurements of the best models for each partition. The quality measurements are the percentage of validation and estimation data unfitness, Akaike's Final Prediction Error (FPE) (Jones, 1975), loss function (Berger, 1985) and mean squared normalized error performance function (MSE) (Poli and Cirillo, 1993). The quality measurements are represented in the graph by 'Val unfit', 'Est unfit', FPE, 'Loss Fcn', and MSE respectively. Although, the model choice in this study is not mathematically based on FPE, loss function and MSE their values were calculated and shown providing an extra measurement of model quality. A variety of measurements is useful for comparing different models as well as comparing the models with different modeling approaches. Differently from the models one to seven, the models for the segments eight and nine present high unfit levels to their validation segments, 72,64% and impossible to calculate, respectively. For $G_8$, the input-output relation of the validation data, shown in the segment eight to nine in Figure 1, is extremely distinct from the data used in the estimation, segment one to eight. For $G_9$, there is only one sample remaining to validate the model. Hence, the model $G_9(s) = 0.6646 / (s + 0.6687)$ corresponds to the structure that better reproduces the sample data with about 93% fitness. Accordingly, $G_9$ is used later to demonstrated the step response.

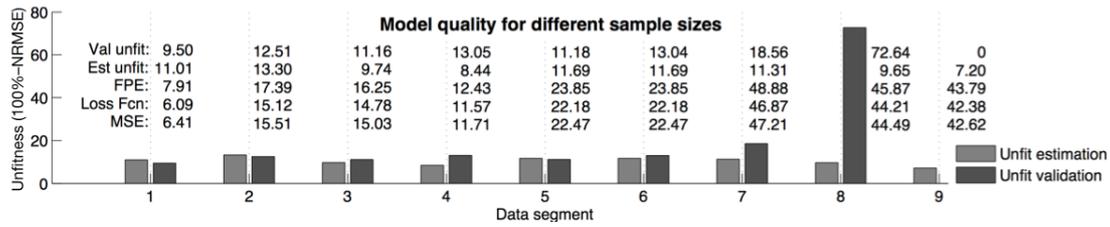

*Figure 2: Quality comparison of the models*

Despite the model $G_9$ has the lowest unfitted data, 'Est unfit', almost the whole sample data was used to estimate $G_9$. Hence, $G_9$ already 'knows' the data sample and for this reason cannot be used as a predictor. In order to illustrate the prediction accuracy of the models, the model with the largest 'unknown' data, *i.e.*, $G_1$ is used. Figure 3 shows the comparison of the measured data and model $G_1(t)$, result of inverse Laplace Transform of $G_1(s) = 0.4193 / (s + 0.4103)$; the solid line is the measured output and the dashed line the model response with 9.5% unfitness. The transfer function $G_1(t)$, can represent the process input-output relationship with sufficient precision. Furthermore, the model is estimated at an early stage, around the initial 10% duration, independently of any previous process knowledge.

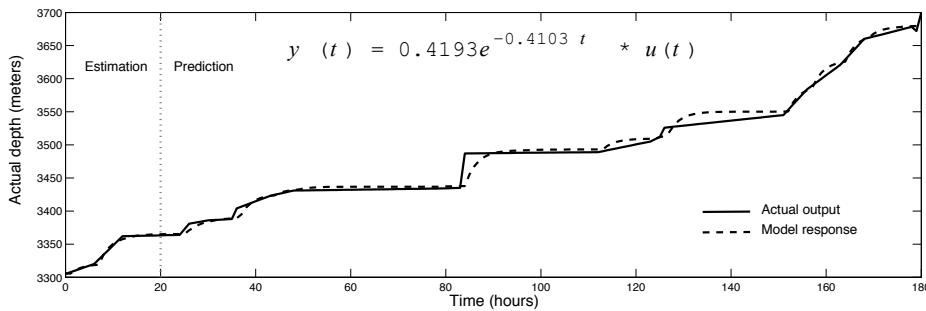

*Figure 3: Comparison between the $G_1$ response and measured data*

## STEP RESPONSE

Figure 4(a) shows the step response for the model $G_9$. The model reaches steady-state about the sixth hour for a threshold of absolute two percent about the final value. The step amplitude used as input, 2.06, is the average drilling goal ahead. The system responds to this input reaching and staying steady at the output peak, $y_p = 2.05$, about the 16th hour, $t_p$. In this case, the steady-state value is the peak value because it is the value that the system tends to when $t \to \infty$ (Mandal, 2006). The average drilling rate from the measured output data is 2.1 meters per hour approximately the model output at steady-state, with a three percent error. Consequently, the model represents the system at steady-state. Although, the transient response stands for a significant part of the system dynamics. The system has a transient response every time it starts or stops. Although it stays at the transient state when it need small corrections, as, for instance, to fit a casing pipe to secure the well. In this case, the system also has inputs as small as one. For this input, the average response of the system is 0.55 meters per hour. In order to assess the system transient response, a unitary step unit was introduced to the system producing the system response, as shown in Figure 4(b). The average response for the unitary input is achieved around one hour by the system indicating that the system is performing in the transient state.

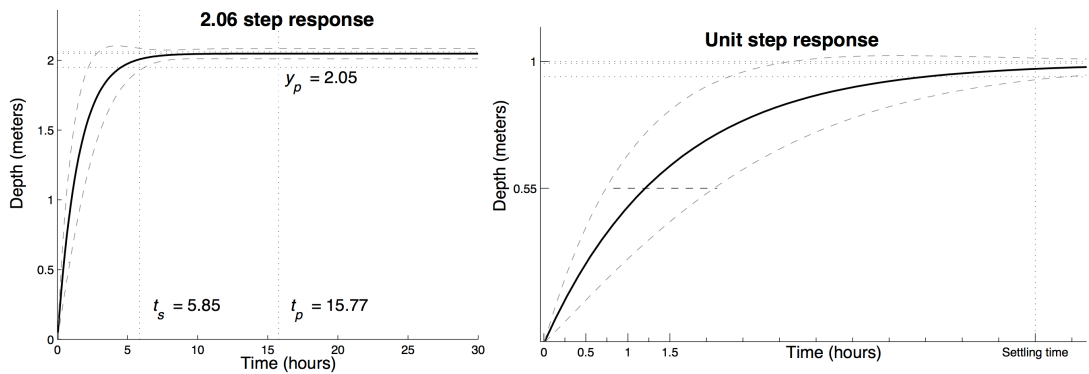

*Figure 4: (a) 2.06 step response (b) Unit step response*

# CONCLUSIONS

The results of the model's accuracy were explicit. The models were consistent with the modeling approach and methodology. The valid transfer functions obtained reliably described the process behavior and presented evidence of their accuracy using a range of model quality measurements. These findings thus lend support to the use of transfer function as a valid model approach and analytical technique in order to describe the dynamic conditions of production in repetitive processes in projects.

Accounting for transient responses, transfer function fulfills a gap left by network scheduling and queueing theory as well as linear and dynamic programming, which ignore the transient stage and assume that the process is at steady-state. Moreover, a transfer functions may act as a multi-level management tool. Because transfer functions provide an output function from an input function, they enable the creation of accurate plans rather than single actions and a throughput function, instead of a system position. Transfer functions may be used by site managers as a process descriptor to monitor and control low-level activities, as shown in Figure 3. Dynamic and accurate plans that respond to actual inputs can regain the trust of those conducting the project work on planning and scheduling. Moreover, the model simulation may be used in a means-ends analysis determining the best solution to a construction process, which frequently requires the optimization of resources to the detriment of shorter duration. In other words, managers may use the model adjusting the drill ahead goal plan until attaining the defined goal, supporting managers' decision-making process. Once the managers are satisfied with both the drill ahead goal plan and the system's outcome, the plan is executed. A transfer function may also be applied to represent higher levels, providing project managers a holistic view.

Reliance on this method must be tempered, however because the case does not represent the general conditions of repetitive process in construction. There is a variety of construction processes that happen in different states, production runs, and different learning curves creating unique process' characteristics. For instance, this study presents the analysis of system's transient and steady-state response, but not unsteady-state because the case scenario does not have this characteristic. Although the model can be reused in similar processes as an initial model, a limitation places on the existence of the process' input-output data. It means that the model accuracy cannot be evaluated until some data has been produced. Finally, the study explores several concepts that are unfamiliar to general construction managers at this point restricting its audience. Nevertheless, the search for and the aggregation of knowledge and expertise from different disciplines and technical fields constitutes the foremost forces driving the evolution in managerial sciences.

**FUTURE OUTCOMES**

Different system identification approaches can write equations for practically any process. However, only after extensive research about the dynamic conditions of production in project-driven systems the lack of knowledge about the transient and unsteady-state responses can be replaced by explanatory and mathematical laws to production in projects. In a further horizon, processes transient and unsteady-state will be understood and managed to generate an optimum process outcome. Being it reducing the transient time, and faster-moving processes to steady-state or applying unsteady-state processing techniques producing an average output above steady-state levels and then creating high-performance processes.